# A new approach to non-lightlike curve pairs

**Zehra Ekinci[1], Mehmet Önder[2]**

*[1]Güzelyurt Locality, 5805 Street, No: 7B/12, 45030, Yunusemre, Manisa, Turkey.*

*E-mail: ari.zehra@windowslive.com*

*[2]Delibekirli Village, Tepe Street, No: 63, 31440, Kırıkhan, Hatay, Turkey.*

*E-mail: mehmetonder197999@gmail.com*

**Abstract**

In this paper, we introduce a new approach to non-lightlike curve pairs by using integral curves in Minkowski 3-space $E_1^3$. We consider direction curve and donor curve to study non-lightlike curve couples such as involute-evolute curves, Mannheim partner curves and Bertrand partner curves. We obtain new methods to construct partner curves of a unit speed non-lightlike curve and give some applications related to helices and slant helices.

## 1. Introduction

Lorentzian geometry plays an important role in the study of special relativity and connected with hyperbolic geometry. It has always preserved its attraction and many scientists from different branches have studied in this area. Famous physicists and mathematicians wrote about the connection between physics and Lorentzian geometry; so did Kant, Helmholtz, Poincaré, Einstein and Hilbert. From the viewpoint of mechanics and cosmology, the geometrical properties of space and time have been considered by Liebscher in Minkowski space[8]. Moreover, when analyzing the characteristics of space curves in Lorentzian geometry, one of the most important issues is the relationships corresponding to the Euclidean case. Then, the well-known special curves in Euclidean geometry such as helix, slant helix, involute-evolute curves, Bertrand curves and Mannheim curves are also fascinating subjects in Minkowski 3-space $E_1^3$. Similar to the Euclidean case, a general helix in $E_1^3$ is a curve for which the ratio of its non-zero curvature to its non-zero torsion is a constant function [3]. Moreover, a slant helix in $E_1^3$ is another special curve which has the property that its principal normal vectors always make a constant angle with a fixed line named its axis. Characterizations of slant helices in $E_1^3$ has been presented by Ali and Lopez by the differential equation of its curvature $\kappa$ and its torsion $\tau$ given by

$$\frac{\kappa^2}{\left(\tau^2-\kappa^2\right)^{3/2}}\left(\frac{\tau}{\kappa}\right)' = \text{constant} \quad \text{or} \quad \frac{\kappa^2}{\left(\kappa^2\pm\tau^2\right)^{3/2}}\left(\frac{\tau}{\kappa}\right)' = \text{constant}$$

([1]). Later, Ali and Turgut have studied position vector of a timelike slant helix [2]. Furthermore, special curve pairs such as involute-evolute curves, Bertrand curves and Mannheim curves are studied largely in Minkowski 3-space [6,7,9].

Moreover, Choi and Kim have introduced a new approach to curve pairs, called principal direction and principal donor curves [4]. Later, Choi, Kim and Ali have studied these donor curves in $E_1^3$ and they have obtained characterizations connected to general helices and slant helices [5]. Önder has defined some new direction curves in Euclidean 3-space and given some relationships between these curves and some special curves such as helix and slant helix [11].

In this study, for non-lightlike curves, we consider integral curves to study partner curves such as involute-evolute curves, Mannheim partner curves and Bertrand partner curves in $E_1^3$. We define these non-lightlike curves as direction curves and obtain relationships between their Frenet elements. Furthermore, we give some applications related to helices and slant helices.

## 2. Preliminaries

The three-dimensional Minkowski space $E_1^3$ is the real vector space $\mathbb{R}^3$ endowed with the standard flat Lorentz metric defined by $\langle , \rangle = -dx_1^2 + dx_2^2 + dx_3^2$ where $(x_1, x_2, x_3)$ is a rectangular coordinate system of $E_1^3$. The Lorentz vector product of two vectors $x = (x_1, x_2, x_3)$, $y = (y_1, y_2, y_3)$ in $E_1^3$ is given by

$$x \times y = (x_2 y_3 - x_3 y_2,\ x_3 y_1 - x_1 y_3,\ x_2 y_1 - x_1 y_2)$$

A vector $x \in E_1^3$ is called a spacelike vector if $\langle x, x \rangle > 0$ or $x = 0$; called a timelike vector if $\langle x, x \rangle < 0$ and called a null (lightlike) vector if $\langle x, x \rangle = 0,\ x \neq 0$. Similarly, a curve $\alpha(t) : I \subset \mathbb{R} \to E_1^3$ is called a spacelike, timelike or null (lightlike) curve if all of its velocity vectors $\alpha'(t)$ are spacelike, timelike or null (lightlike), respectively [10]. For a non-lightlike curve $\alpha$, the parameter $s$ defined by $s = \int_0^s \|\alpha'(t)\| dt$ is called arc-length parameter and if $\|\alpha'(t)\| = 1$, then $\alpha$ is called unit speed curve, where $\|\alpha'(t)\| = \sqrt{|\langle \alpha'(t), \alpha'(t) \rangle|}$.

Let $\{T, N, B\}$ denotes the Frenet frame of a non-lightlike curve $\alpha$. Then, the curve $\alpha$ is a timelike or a spacelike curve. The spacelike curves with non-lightlike Frenet vectors have two types according to the Lorentzian casual characters of Frenet vectors: A spacelike curve is called of type 1 (respectively, type 2) if its principal normal vector $N$ (respectively, binormal vector $B$) is timelike and other Frenet vectors are spacelike [5]. For the derivatives of the Frenet frame of a unit speed non-lightlike curve $\alpha(s)$, the following Frenet-Serret formulae hold:

$$\begin{bmatrix} T' \\ N' \\ B' \end{bmatrix} = \begin{bmatrix} 0 & \kappa & 0 \\ \varepsilon_B \kappa & 0 & \tau \\ 0 & \varepsilon_T \tau & 0 \end{bmatrix} \begin{bmatrix} T \\ N \\ B \end{bmatrix} \tag{1}$$

where $\varepsilon_Y = \langle Y, Y \rangle$, $\varepsilon_Y^2 = 1$, $Y' = dY/ds$, $Y \in \{T, N, B\}$, $\kappa(s)$ is the curvature and $\tau(s)$ is the torsion (or second curvature) of $\alpha$ at $s$ defined by $\tau(s) = \varepsilon_N \varepsilon_B \langle B', N \rangle$. The curve $\alpha(s)$ is called Frenet curve if $\kappa(s) \neq 0$ [12]. The Frenet vectors satisfy the relations

$$B = \varepsilon_T \varepsilon_N T \times N,\quad N = \varepsilon_B \varepsilon_T B \times T,\quad T = \varepsilon_N \varepsilon_B N \times B,\quad \varepsilon_B = -\varepsilon_T \varepsilon_N .$$

Now, we give the definitions of some associated curves defined by Choi and Kim [5]. Let $I \subset \mathbb{R}$ be an open interval. For a non-lightlike Frenet curve $\alpha : I \to E_1^3$, consider a non-lightlike vector field $X$ given by

$$X(s) = u(s)T(s) + v(s)N(s) + w(s)B(s),$$

where $u, v$ and $w$ are arbitrary differentiable functions of $s$ which is the arc-length parameter of $\alpha$. In physics, if $X$ is a force field, the integral curves of $X$ are called lines of force. If $X$ is the velocity of the fluid flow, the integral curves of $X$ are called lines of flow. These are the paths of motion of the fluid particles.

Let $\varepsilon_T u^2(s)+\varepsilon_N v^2(s)+\varepsilon_B w^2(s)=\sigma(s)=\pm 1$ holds. Then the definitions of $X$-direction curve and $X$-donor curve in $E_1^3$ are given as follows.

**Definition 2.2. ([5])** Let $\alpha$ be a non-lightlike Frenet curve in $E_1^3$ with Frenet frame $\{T,N,B\}$ and $X$ be a unit vector field along $\alpha$. The curve $\beta: I\to E_1^3$ is called a $X$-direction curve of $\alpha$, if the tangent $\overline{T}$ of $\beta$ is equal to $X$, i.e., $\beta$ is an integral curve of $X$. The curve $\alpha$ whose $X$-direction curve is $\beta$ is called the $X$-donor curve of $\beta$ in $E_1^3$.

Let $\alpha: I\to E_1^3$ be a unit speed non-lightlike Frenet curve with arc-length parameter $s$, Frenet frame $\{T,N,B\}$ and curvatures $\kappa,\ \tau$. And let $X$ be a non-lightlike continuous vector valued function along $\alpha$ defined by

$$X(s)=u(s)T(s)+v(s)N(s)+w(s)B(s), \tag{2}$$

where $u(s),\ v(s),\ w(s)$ are differentiable functions on $I$ satisfying

$$\langle X,X\rangle=\varepsilon_T u^2(s)+\varepsilon_N v^2(s)+\varepsilon_B w^2(s)=\sigma(s)=\pm 1, \tag{3}$$

and $\beta: I\to E_1^3$ be an $X$-direction curve of $\alpha$ in $E_1^3$. The Frenet vectors and curvatures of $\beta$ be denoted by $\{\overline{T},\overline{N},\overline{B}\}$ and $\overline{\kappa},\ \overline{\tau}$, respectively. From (2) and (3), it is clear that the arc-length parameter $s$ of $\alpha$ can be also taken as arc-length parameter of $\beta$. Then, differentiating (2) with respect to $s$ and using the fact that $X=\overline{T}$, it follows

$$\overline{\kappa}\overline{N}=(u'+v\varepsilon_B\kappa)T+(v'+u\kappa+w\varepsilon_T\tau)N+(w'+v\tau)B \tag{4}$$

Now, by using equality (4), we study non-lightlike curve pairs such as involute-evolute curves, Mannheim partner curves and Bertrand partner curves. In the following sections, we will assume that $X$ is a non-lightlike continuous vector valued function along $\alpha$ as given in (2) and it satisfies (3).

## 3. Non-lightlike involute-evolute-direction curves

In this section, we will give definitions of non-lightlike involute-evolute-direction curves and obtain relationships between these curves.

**Definition 3.1.** Let $\alpha: I\to E_1^3$ be a unit speed non-lightlike Frenet curve and $\beta: I\to E_1^3$ be an $X$-direction curve of $\alpha$. If $\beta$ is an evolute of $\alpha$ in $E_1^3$, then $\beta$ is called evolute-direction curve of $\alpha$ and $\alpha$ is said to be involute-donor curve of $\beta$.

**Theorem 3.1.** *For the non-lightlike Frenet curve $\alpha$, the curve $\beta$ is an evolute-direction curve of $\alpha$ if and only if followings hold,*

$$\text{(i)}\ u(s)=0,\ v(s)=-\cos\left(\int\tau ds\right),\ w(s)=\sin\left(\int\tau ds\right), \tag{5}$$

*if $\alpha$ is a timelike curve and $X$ is a spacelike vector field.*

$$\text{(ii)}\ u(s)=0,\ v(s)=-\cosh\left(\int\tau ds\right),\ w(s)=\sinh\left(\int\tau ds\right), \tag{6}$$

*if $\alpha$ is a spacelike curve of type 1 (respectively, type 2) and $X$ is timelike (respectively, spacelike).*

$$\text{(iii)}\ u(s)=0,\ v(s)=\sinh\left(\int\tau ds\right),\ w(s)=-\cosh\left(\int\tau ds\right), \tag{7}$$

*if $\alpha$ is a spacelike curve of type 1 (respectively, type 2) and $X$ is spacelike (respectively timelike)*

**Proof.** Let $\beta$ be an evolute-direction curve of $\alpha$. From the definition of involute-evolute curves, it is clear that $\langle \bar{T}, T \rangle = 0$ and $\bar{N} = T$. Then, from (4) we have that $\beta$ is non-lightlike evolute-direction curve of $\alpha$ if and only if $u = 0$ and from (3) and (4) we can write the system

$$\varepsilon_B v\,\kappa = \bar{\kappa} \neq 0,\quad v' + \varepsilon_T w\,\tau = 0,\quad w' + v\tau = 0\,;\ \ \varepsilon_N v^2(s) + \varepsilon_B w^2(s) = \sigma(s)\,. \tag{8}$$

In order to solve this system, we must consider the Lorentzian casual characters of involute-evolute-direction curves and $X$ as follows:

(i) If $\alpha$ is a timelike curve and $X$ is spacelike, from (8), it follows

$$v' - w\tau = 0,\ \ w' + v\tau = 0,\ \ v^2(s) + w^2(s) = 1\,. \tag{9}$$

The solution of the above system is

$$\left\{u(s) = 0,\ v(s) = -\cos\left(\int \tau ds\right),\ w(s) = \sin\left(\int \tau ds\right)\right\}. \tag{10}$$

(ii) If $\alpha$ is a spacelike curve of type 1 (respectively, type 2) and $X$ is timelike (respectively, spacelike), then equations (8) become

$$v' + w\tau = 0,\ \ w' + v\tau = 0,\ \ v^2(s) - w^2(s) = 1 \tag{11}$$

and the solution is

$$\left\{u(s) = 0,\ v(s) = -\cosh\left(\int \tau ds\right),\ w(s) = \sinh\left(\int \tau ds\right)\right\}. \tag{12}$$

(iii) If $\alpha$ is a spacelike curve of type 1 (respectively, type 2) and $X$ is spacelike (respectively, timelike), then from equations (8), we have

$$v' + w\tau = 0,\ \ w' + v\tau = 0,\ \ -v^2(s) + w^2(s) = 1\,, \tag{13}$$

and the solution is obtained as

$$\left\{u(s) = 0,\ v(s) = \sinh\left(\int \tau ds\right),\ w(s) = -\cosh\left(\int \tau ds\right)\right\}. \tag{14}$$

Theorem 3.1 allows us to give the following definition.

**Definition 3.2.** Let $\alpha$ be a non-lightlike curve. An integral curve of one of the vector fields

$$-\cos\left(\int \tau ds\right) N(s) + \sin\left(\int \tau ds\right) B(s),\quad -\cosh\left(\int \tau ds\right) N(s) + \sinh\left(\int \tau ds\right) B(s)$$

$$\sinh\left(\int \tau ds\right) N(s) - \cosh\left(\int \tau ds\right) B(s)\,,$$

is called non-lightlike evolute-direction curve of $\alpha$.

From Theorem 3.1, we obtain a method to construct an evolute of a given non-lightlike curve by using its Frenet vectors $N, B$ and its torsion $\tau$. It means that to obtain a non-lightlike evolute-direction curve of a non-lightlike curve, it is enough to know the Frenet elements $N, B, \tau$ of reference curve.

Moreover, From Theorem 3.1, the relationships between curvatures and Frenet vectors of non-lightlike involute-evolute-direction curves can be given by the following corollary.

**Corollary 3.1.** *If $\alpha$ is a timelike curve and $X$ is spacelike, then the relations between the Frenet vectors of involute-evolute-direction curves are given as follows:*

$$\bar{T}(s) = -\cos\left(\int \tau ds\right) N(s) + \sin\left(\int \tau ds\right) B(s),\qquad \bar{N}(s) = T(s), \tag{15}$$

$$\bar{B}(s) = -\sin\left(\int \tau ds\right) N(s) - \cos\left(\int \tau ds\right) B(s)\,. \tag{16}$$

*For the other cases given in Theorem 3.1, the relations between the Frenet vectors can be*

*given easily.*

**Theorem 3.2.** *Let $\alpha$ be a non-lightlike Frenet curve and $\beta$ be an evolute-direction of $\alpha$.*

*(i) If $\alpha$ is a timelike curve, $X$ is spacelike and $\beta$ is a spacelike evolute-direction curve of type 1, then*

$$\overline{\kappa} = \kappa\left|\cos\left(\int \tau ds\right)\right|,\ \overline{\tau} = -\kappa\sin\left(\int \tau ds\right),\ \ \kappa = \sqrt{\overline{\kappa}^2 + \overline{\tau}^2}\ ,\ \tau = \frac{\overline{\kappa}^2}{\overline{\kappa}^2 + \overline{\tau}^2}\left(\frac{\overline{\tau}}{\overline{\kappa}}\right)'. \qquad (17)$$

*(ii) If $\alpha$ is a spacelike curve of type 1 (respectively, type 2), $X$ is timelike (respectively, spacelike) and $\beta$ is a timelike evolute-direction curve with $|\overline{\kappa}| > |\overline{\tau}|$ (respectively, spacelike curve of type 2), then*

$$\overline{\kappa} = \kappa\left|\cosh\left(\int \tau ds\right)\right|,\ \overline{\tau} = -\kappa\sinh\left(\int \tau ds\right),\ \ \kappa = \sqrt{\overline{\kappa}^2 - \overline{\tau}^2}\ ,\ \tau = \frac{\overline{\kappa}^2}{\overline{\kappa}^2 - \overline{\tau}^2}\left(\frac{\overline{\tau}}{\overline{\kappa}}\right)' \qquad (18)$$

*(iii) If $\alpha$ is a spacelike curve of type 1 (respectively, type 2), $X$ is spacelike (respectively, timelike) and $\beta$ is a spacelike evolute-direction curve of type 2 with $|\overline{\tau}| > |\overline{\kappa}|$, (respectively, timelike) then*

$$\overline{\kappa} = \kappa\left|\sinh\left(\int \tau ds\right)\right|,\ \overline{\tau} = \kappa\cosh\left(\int \tau ds\right),\ \ \kappa = \sqrt{\overline{\tau}^2 - \overline{\kappa}^2}\ ,\ \ \tau = \frac{\overline{\kappa}^2}{\overline{\tau}^2 - \overline{\kappa}^2}\left(\frac{\overline{\tau}}{\overline{\kappa}}\right)' \qquad (19)$$

**Proof.** Now, we give the proof of case (i). The proofs of (ii) and (iii) can be given by a similar way. From (5) and the first equation of system (8), we have

$$\overline{\kappa} = \kappa\left|\cos\left(\int \tau ds\right)\right|. \qquad (20)$$

Moreover from (16), we easily get $\overline{B}' = -\kappa\sin\left(\int \tau ds\right)T$. Since $\overline{N} = T$ and $\alpha$ is timelike, we have

$$\overline{\tau} = \left\langle \overline{B}', \overline{N} \right\rangle = -\kappa\sin\left(\int \tau ds\right). \qquad (21)$$

From (20) and (21), it follows $\kappa = \sqrt{\overline{\kappa}^2 + \overline{\tau}^2}$. Substituting that in (20) and (21), we have

$$\sin\left(\int \tau ds\right) = -\frac{\overline{\tau}}{\sqrt{\overline{\kappa}^2 + \overline{\tau}^2}},\quad \cos\left(\int \tau ds\right) = -\frac{\overline{\kappa}}{\sqrt{\overline{\kappa}^2 + \overline{\tau}^2}}, \qquad (22)$$

respectively. Differentiating first equality in (22) with respect to $s$, we have

$$\tau\cos\left(\int \tau ds\right) = \frac{\overline{\kappa}(\overline{\kappa}'\overline{\tau} - \overline{\kappa}\,\overline{\tau}')}{(\overline{\kappa}^2 + \overline{\tau}^2)^{3/2}}. \qquad (23)$$

From (20) and second equality in (22), it follows

$$\tau = \frac{\overline{\kappa}'\,\overline{\tau} - \overline{\kappa}\,\overline{\tau}'}{\overline{\kappa}^2 + \overline{\tau}^2}, \qquad (24)$$

or equivalently,

$$\tau = \frac{\overline{\kappa}^2}{\overline{\kappa}^2 + \overline{\tau}^2}\left(\frac{\overline{\tau}}{\overline{\kappa}}\right)'. \qquad (25)$$

**Corollary 3.2.** *Let $\alpha$ be a non-lightlike Frenet curve and $\beta$ be a non-lightlike evolute-direction curve of $\alpha$. Then, one of the following equalities hold,*

$$(i)\ \frac{\tau}{\kappa}=\frac{\bar{\kappa}^2}{\left(\bar{\kappa}^2+\bar{\tau}^2\right)^{3/2}}\left(\frac{\bar{\tau}}{\bar{\kappa}}\right)' \quad (ii)\ \frac{\tau}{\kappa}=\frac{\bar{\kappa}^2}{\left(\bar{\tau}^2-\bar{\kappa}^2\right)^{3/2}}\left(\frac{\bar{\tau}}{\bar{\kappa}}\right)' \quad (iii)\ \frac{\tau}{\kappa}=\frac{\bar{\kappa}^2}{\left(\bar{\kappa}^2-\bar{\tau}^2\right)^{3/2}}\left(\frac{\bar{\tau}}{\bar{\kappa}}\right)' \qquad (26)$$

Considering Corollary 3.2, we give the following theorems which give two ways to construct special curves: 1) constructing non-lightlike slant helices by using non-lightlike general helices; 2) constructing non-lightlike helices by using non-lightlike plane curves.

***Theorem 3.3.*** *Let $\alpha : I \to E_1^3$ be a non-lightlike curve and $\beta$ be a non-lightlike evolute-direction curve of $\alpha$. Then the followings are equivalent,*

*(i) $\alpha$ is a helix curve.*

*(ii) $\alpha$ is an involute-donor curve of a slant helix.*

*(iii) An evolute-direction curve of $\alpha$ is a slant helix.*

***Theorem 3.4.*** *Let $\alpha : I \to E_1^3$ be a non-lightlike curve and $\beta$ be a non-lightlike evolute-direction curve of $\alpha$. The following statements are equivalent,*

*(i) $\alpha$ is a plane curve.*

*(ii) $\alpha$ is an involute-donor curve of a helix.*

*(iii) An evolute-direction curve of $\alpha$ is a helix.*

**4. Non-lightlike Mannheim-direction curves**

In this section, we define non-lightlike Mannheim-direction curves and obtain relationships between these curves in $E_1^3$.

**Definition 4.1.** Let $\alpha : I \to E_1^3$ be a non-lightlike Frenet curve and $\beta : I \to E_1^3$ be non-lightlike $X$-direction curve of $\alpha$. If $\beta$ is a Mannheim curve of $\alpha$ and $\alpha$ is a Mannheim partner curve of $\beta$, then $\beta$ is called Mannheim-direction curve of $\alpha$ and $\alpha$ is said to be Mannheim-donor curve of $\beta$.

**Theorem 4.1.** *For the Frenet curve $\alpha : I \to E_1^3$, the curve $\beta : I \to E_1^3$ is a Mannheim-direction curve of $\alpha$ if and only if*

$$(i)\ u(s)=-\cosh\left(\int \kappa ds\right),\ v(s)=\sinh\left(\int \kappa ds\right),\ w(s)=0, \qquad (27)$$

*if $\alpha$ is a timelike curve (respectively, spacelike curve of type 1) and $X$ is timelike (respectively, spacelike).*

$$(ii)\ u(s)=\sinh\left(\int \kappa ds\right),\ v(s)=-\cosh\left(\int \kappa ds\right),\ w(s)=0, \qquad (28)$$

*if $\alpha$ is a timelike curve (respectively, spacelike curve of type 1) and $X$ is spacelike (respectively, timelike).*

$$(iii)\ u(s)=-\cos\left(\int \kappa ds\right),\ v(s)=\sin\left(\int \kappa ds\right),\ w(s)=0, \qquad (29)$$

*if $\alpha$ is a spacelike curve of type 2 and $X$ is spacelike.*

**Proof.** Let $\beta$ be a Mannheim curve of $\alpha$. From the definition of Mannheim curves, it is well-known that $\bar{N}=B$ [7,9]. Then, from (4), we have that $\beta$ is a Mannheim-direction curve of $\alpha$ if and only if the system

$$u'+v\varepsilon_B\kappa=0,\quad u\kappa+v'+w\varepsilon_T\tau=0,\quad w'+v\tau=\bar{\kappa}\neq 0, \qquad (30)$$

holds. Multiplying the first equation in (30) with $\varepsilon_T u$, second equation with $\varepsilon_N v$ and adding the results give

$$\varepsilon_T u u' + \varepsilon_N v v' + w \varepsilon_N \varepsilon_T \tau = 0 \tag{31}$$

From (3) and (31), we have $\varepsilon_B w(w' + v\tau) = 0$. Since $w' + v\tau = \overline{\kappa} \neq 0$ in (30) and $\varepsilon_B \neq 0$, it follows $w = 0$. Then the system (30) reduced to the following system

$$u' + v\,\varepsilon_B\,\kappa = 0, \quad u\kappa + v' = 0, \quad v\tau = \overline{\kappa} \neq 0, \tag{32}$$

and from (3)

$$\varepsilon_T u^2(s) + \varepsilon_N v^2(s) = \sigma(s) = \pm 1. \tag{33}$$

In order to solve the system, we consider the Lorentzian casual characters of curves and continuous vector valued function $X$ as follows:

(i) If $\alpha$ is a timelike curve (respectively, spacelike curve of type 1) and $X$ is a timelike (respectively, spacelike), then (32) and (33) give

$$u' + v\kappa = 0, \quad u\kappa + v' = 0, \quad u^2(s) - v^2(s) = 1. \tag{34}$$

So the solution of the system is

$$\left\{ u(s) = \cosh\left(\int \kappa ds\right), \ v(s) = -\sinh\left(\int \kappa ds\right), \ w(s) = 0 \right\}. \tag{35}$$

(ii) If $\alpha$ is a timelike curve (respectively, spacelike curve of type 1) and $X$ is a spacelike (respectively, timelike), then from (32) and (33) we have

$$u' + v\kappa = 0, \quad u\kappa + v' = 0, \quad -u^2(s) + v^2(s) = 1, \tag{36}$$

which has the solution

$$\left\{ u(s) = \sinh\left(\int \kappa ds\right), \ v(s) = -\cosh\left(\int \kappa ds\right), \ w(s) = 0 \right\}. \tag{37}$$

(iii) If $\alpha$ is a spacelike curve of type 2 and $X$ is spacelike, then from (32) and (33) it follows

$$u' - v\kappa = 0, \quad u\kappa + v' = 0, \quad u^2(s) + v^2(s) = 1, \tag{38}$$

Then, the solution is

$$\left\{ u(s) = -\cos\left(\int \kappa ds\right), \ v(s) = \sin\left(\int \kappa ds\right), \ w(s) = 0 \right\}. \tag{39}$$

**Definition 4.2.** Let $\alpha$ be a non-lightlike curve. An integral curve of one of the vector fields is called Mannheim-direction curve of $\alpha$,

$$\begin{cases} -\cosh\left(\int \kappa ds\right) T(s) + \sinh\left(\int \kappa ds\right) N(s), \quad \sinh\left(\int \kappa ds\right) T(s) - \cosh\left(\int \kappa ds\right) N(s), \\ -\cos\left(\int \kappa ds\right) T(s) + \sin\left(\int \kappa ds\right) N(s) \end{cases} \tag{40}$$

From Theorem 4.1, we obtain a method to construct a non-lightlike Mannheim curve of a unit speed non-lightlike curve by using its Frenet elements $T$, $N$ and $\kappa$.

Now we can give the relationships between curvatures and Frenet vectors of Mannheim-direction curves as follows:

**Corollary 4.1.** *If $\alpha$ is a timelike curve and $X$ is timelike, then the relations between the Frenet vectors of Mannheim-direction curves are given by*

$$\overline{T}(s) = -\cosh\left(\int \kappa ds\right) T(s) + \sinh\left(\int \kappa ds\right) N(s), \quad \overline{N}(s) = B(s), \tag{41}$$

$$\overline{B}(s) = \sinh\left(\int \kappa ds\right) T(s) - \cosh\left(\int \kappa ds\right) N(s). \tag{42}$$

**Proof.** The proof is clear from (27), (28), (29) and definition of non-lightlike Mannheim curves. For the cases (ii), (iii), (iv) and (v) given in Theorem 4.1, the relations between Frenet vectors can be obtained easily.

**Theorem 4.2.** *Let $\alpha$ be a non-lightlike curve and $\beta$ be a Mannheim-direction curve of $\alpha$.*

*(i) If $\alpha$ is a timelike curve with $|\tau|>|\kappa|$ (respectively, spacelike curve of type 1), $X$ is timelike (respectively, spacelike) and $\beta$ is a timelike Mannheim-direction curve (respectively, spacelike curve of type 2) of $\alpha$, then*

$$\bar{\kappa}=\left|\tau\sinh\left(\int\kappa ds\right)\right|,\ \bar{\tau}=-\tau\cosh\left(\int\kappa ds\right),\ \kappa=\left|\frac{\bar{\kappa}^2}{\bar{\tau}^2-\bar{\kappa}^2}\left(\frac{\bar{\tau}}{\bar{\kappa}}\right)'\right|,\ |\tau|=\sqrt{\bar{\tau}^2-\bar{\kappa}^2}\,. \tag{43}$$

*(ii) If $\alpha$ is a timelike curve with $|\kappa|>|\tau|$ (respectively, spacelike curve of type 1), $X$ is spacelike (respectively, timelike) and $\beta$ is a Mannheim-direction spacelike curve of type 2 (respectively, timelike) of $\alpha$, then*

$$\bar{\kappa}=\left|\tau\cosh\left(\int\kappa ds\right)\right|,\ \bar{\tau}=\tau\sinh\left(\int\kappa ds\right),\ \kappa=\left|\frac{\bar{\kappa}^2}{\bar{\kappa}^2-\bar{\tau}^2}\left(\frac{\bar{\tau}}{\bar{\kappa}}\right)'\right|,\ |\tau|=\sqrt{\bar{\kappa}^2-\bar{\tau}^2}\,. \tag{44}$$

*(iii) If $\alpha$ is a spacelike curve of type 2, $X$ is spacelike and $\beta$ is a Mannheim-direction spacelike curve of type 1 of $\alpha$, then*

$$\bar{\kappa}=\left|\tau\sin\left(\int\kappa ds\right)\right|,\ \bar{\tau}=\tau\cos\left(\int\kappa ds\right),\ \kappa=\left|\frac{\bar{\kappa}^2}{\bar{\kappa}^2+\bar{\tau}^2}\left(\frac{\bar{\tau}}{\bar{\kappa}}\right)'\right|,\ |\tau|=\sqrt{\bar{\kappa}^2+\bar{\tau}^2}\,. \tag{45}$$

**Proof.** Now, we give the proof of first case (i). The proofs of the cases (ii) and (iii) can be given by a similar way. From (34) and the first equation of system (39), we have

$$\bar{\kappa}=\left|\tau\sinh\left(\int\kappa ds\right)\right|. \tag{46}$$

Moreover from (42), we easily get $\bar{B}'=-\tau\cosh\left(\int\kappa ds\right)B$. Since $\bar{N}=B$ and $\alpha$ is timelike, we have

$$\bar{\tau}=\left\langle \bar{B}',\bar{N}\right\rangle=-\tau\cosh\left(\int\kappa ds\right). \tag{47}$$

From (46) and (47), we get easily $|\tau|=\sqrt{\bar{\tau}^2-\bar{\kappa}^2}$. Moreover, substituting $|\tau|=\sqrt{\bar{\tau}^2-\bar{\kappa}^2}$ into (46) and (47), it follows

$$\sinh\left(\int\tau ds\right)=\frac{\bar{\kappa}}{\sqrt{\bar{\tau}^2-\bar{\kappa}^2}},\quad \cosh\left(\int\kappa ds\right)=-\frac{\bar{\tau}}{\sqrt{\bar{\tau}^2-\bar{\kappa}^2}}, \tag{48}$$

respectively. Differentiating second equality in (48) with respect to $s$, we have

$$\kappa\sinh\left(\int\kappa ds\right)=\frac{\bar{\kappa}(\bar{\kappa}\bar{\tau}'-\bar{\kappa}'\bar{\tau})}{(\bar{\tau}^2-\bar{\kappa}^2)^{3/2}}. \tag{49}$$

From (46) and (49), it follows

$$\kappa=\frac{\bar{\kappa}\bar{\tau}'-\bar{\kappa}'\bar{\tau}}{\bar{\tau}^2-\bar{\kappa}^2}, \tag{50}$$

or equivalently,

$$\kappa=\left|\frac{\bar{\kappa}^2}{\bar{\tau}^2-\bar{\kappa}^2}\left(\frac{\bar{\tau}}{\bar{\kappa}}\right)'\right|. \tag{51}$$

**Corollary 4.2.** *Let $\beta$ be Mannheim-direction curve of $\alpha$. Then one of the followings hold*

$$\frac{\tau}{\kappa}=\mp\frac{1}{\left|\frac{\bar{\kappa}^2}{\left(\bar{\tau}^2-\bar{\kappa}^2\right)^{3/2}}\left(\frac{\bar{\tau}}{\bar{\kappa}}\right)'\right|},\quad \frac{\tau}{\kappa}=\mp\frac{1}{\left|\frac{\bar{\kappa}^2}{\left(\bar{\kappa}^2-\bar{\tau}^2\right)^{3/2}}\left(\frac{\bar{\tau}}{\bar{\kappa}}\right)'\right|},\quad \frac{\tau}{\kappa}=\mp\frac{1}{\left|\frac{\bar{\kappa}^2}{\left(\bar{\kappa}^2+\bar{\tau}^2\right)^{3/2}}\left(\frac{\bar{\tau}}{\bar{\kappa}}\right)'\right|}.$$

By using Corollary 4.2, we obtain a method to construct the slant helices by using general helices for which the obtained slant helix is also a Mannheim curve of the reference helix in $E_1^3$. This fact can be given by the following theorem.

***Theorem 4.3.*** *Let $\alpha: I \to E_1^3$ be a non-lightlike curve and $\beta$ be a Mannheim-direction curve of $\alpha$. The following statements are equivalent*

*(i) $\alpha$ is a helix.*

*(ii) $\alpha$ is a Mannheim-donor curve of a slant helix.*

*(iii) A Mannheim-direction curve of $\alpha$ is a slant helix.*

**5. Non-lightlike Bertrand-direction curves**

In this section, we define non-lightlike Bertrand-direction curves and obtain relationships between these curves in $E_1^3$.

**Definition 5.1.** Let $\alpha: I \to E_1^3$ be a non-lightlike Frenet curve and $\beta: I \to E_1^3$ be non-lightlike $X$-direction curve of $\alpha$. If $\beta$ is a Bertrand curve of $\alpha$ and $\alpha$ is a Bertrand partner curve of $\beta$, then $\beta$ is called Bertrand-direction curve of $\alpha$ and is said to be Bertrand-donor curve of $\beta$ in $E_1^3$.

**Theorem 5.1.** *For the non-lightlike Frenet curve $\alpha: I \to E_1^3$, the curve $\beta: I \to E_1^3$ is a Bertrand-direction curve of $\alpha$ if and only if*

*(i)* $u(s)=\cosh\theta,\ v(s)=0,\ w(s)=\sinh\theta$, (52)

*if $\alpha$ is a timelike curve (respectively, spacelike curve of type 2) and $X$ is timelike (respectively, spacelike).*

*(ii)* $u(s)=\sinh\theta,\ v(s)=0,\ w(s)=\cosh\theta$, (53)

*if $\alpha$ is a timelike curve (respectively, spacelike curve of type 2) and $X$ is spacelike (respectively, timelike).*

*(iii)* $u(s)=\cos\theta,\ v(s)=0,\ w(s)=\sin\theta$, (54)

*if $\alpha$ is a spacelike curve of type 1 and $X$ is spacelike. In all cases $\theta$ is the constant angle between the tangent lines of the curves.*

**Proof.** From the definition of Bertrand curves, it is well-known that $\bar{N}=N$. Then, from (4) we have that $\beta$ is a Bertrand-direction curve of $\alpha$ in $E_1^3$ if and only if the system,

$$u'+v\varepsilon_B\kappa=0,\quad u\kappa+v'+w\varepsilon_T\tau=\bar{\kappa}\neq 0,\quad w'+v\tau=0, \tag{55}$$

holds. Multiplying the first equation in (55) with $\varepsilon_T u$ and second equation with $\varepsilon_N v$ and third equation with $\varepsilon_B w$ and adding the results gives

$$\varepsilon_T uu'+\varepsilon_N vv'+\varepsilon_B ww'=\varepsilon_N v\bar{\kappa}. \tag{56}$$

From (3) and (56), we have $\varepsilon_N v\bar{\kappa}=0$. Since $\bar{\kappa}\neq 0$ and $\varepsilon_N\neq 0$, it follows $v=0$. Then the system (55) reduced to the following system

$$u'=0,\quad u\kappa+w\varepsilon_T\tau=\bar{\kappa}\neq 0,\quad w'=0, \tag{57}$$

and from (3)

$$\varepsilon_T u^2(s) + \varepsilon_B w^2(s) = \sigma(s)\,. \tag{58}$$

In order to solve this system, we should consider the Lorentzian casual characters of Bertrand-direction curves and continuous vector valued function $X$ as follows:

(i) If $\alpha$ is a timelike curve (respectively, spacelike curve of type 2) and $X$ is timelike (respectively, spacelike) then (57) and (58) give

$$u' = 0, \quad w' = 0, \quad u^2(s) - w^2(s) = 1 \tag{59}$$

So the solution of the above system is

$$\{u(s) = \cosh\theta,\ v(s) = 0,\ w(s) = \sinh\theta;\ \ \theta = constant\}\,. \tag{60}$$

(ii) If $\alpha$ is a timelike curve (respectively, spacelike curve of type 2) and $X$ is spacelike (respectively, timelike) then (57) and (58) become

$$u' = 0, \quad w' = 0, \quad -u^2(s) + w^2(s) = 1 \tag{61}$$

So the solution of the above system is

$$\{u(s) = \sinh\theta,\ v(s) = 0,\ w(s) = \cosh\theta;\ \ \theta = constant\}\,. \tag{62}$$

(iii) If $\alpha$ is a spacelike curve of type 1 and $X$ is spacelike then from (57) and (58) it follows

$$u' = 0, \quad w' = 0, \quad u^2(s) + w^2(s) = 1\,. \tag{63}$$

So the solution is

$$\{u(s) = \cos\theta,\ v(s) = 0,\ w(s) = \sin\theta;\ \ \theta = constant\}\,. \tag{64}$$

**Definition 5.2.** Let $\alpha$ be a non-lightlike Frenet curve. An integral curve of one of the vector fields

$$\cosh\theta T(s) + \sinh\theta B(s)\,, \quad \sinh\theta T(s) + \cosh\theta B(s)\,, \quad \cos\theta T(s) + \sin\theta B(s) \tag{65}$$

is called Bertrand-direction curve of $\alpha$, where $\theta = \text{constant}$.

From Theorem 5.1 we obtain a method to construct a non-lightlike Bertrand curve of a unit speed non-lightlike curve by using its Frenet vectors $T, B$ and a constant $\theta$ in $E_1^3$.

The relationships between curvatures and Frenet vectors of non-lightlike Bertrand-direction curves in $E_1^3$ can be given as follows.

**Corollary 5.1.** *If $\alpha$ is a timelike curve and $X$ is timelike, then*

$$\overline{T}(s) = \cosh\theta T(s) + \sinh\theta B(s), \quad \overline{N}(s) = N(s), \tag{66}$$

$$\overline{B}(s) = -\sinh\theta T(s) - \cosh\theta B(s)\,. \tag{67}$$

**Proof.** The proof is clear from Theorem 5.1. For other cases given in Theorem 5.1, corresponding equalities for (66) and (67) can be obtained easily.

**Theorem 5.2.** *Let $\alpha$ be a non-lightlike Frenet curve and $\beta$ be a non-lightlike Bertrand-direction curve of $\alpha$.*

*(i) If $\alpha$ is a timelike curve, $X$ is timelike and $\beta$ is a timelike Bertrand-direction curve of $\alpha$, then*

$$\overline{\kappa} = \left|\kappa\cosh\theta - \tau\sinh\theta\right|, \quad \overline{\tau} = -\kappa\sinh\theta + \tau\cosh\theta\,. \tag{68}$$

$$\kappa = \left|\overline{\kappa}\cosh\theta + \overline{\tau}\sinh\theta\right|, \ \tau = \overline{\kappa}\sinh\theta + \overline{\tau}\cosh\theta\,. \tag{69}$$

*(ii) If $\alpha$ is a timelike curve, $X$ is spacelike and $\beta$ is a spacelike Bertrand-direction curve of type 2 of $\alpha$, then*

$$\bar{\kappa} = \left|\kappa \sinh\theta - \tau \cosh\theta\right|, \quad \bar{\tau} = -\kappa \cosh\theta + \tau \sinh\theta. \tag{70}$$

$$\kappa = \left|-\bar{\kappa} \sinh\theta - \bar{\tau} \cosh\theta\right|, \tau = -\bar{\kappa} \cosh\theta - \bar{\tau} \sinh\theta. \tag{71}$$

*(iii) If $\alpha$ is a spacelike curve of type 1, $X$ is spacelike and $\beta$ is a spacelike Bertrand-direction curve of type 1, then*

$$\bar{\kappa} = \left|\kappa \cos\theta + \tau \sin\theta\right|, \quad \bar{\tau} = -\kappa \sin\theta + \tau \cos\theta. \tag{72}$$

$$\kappa = \left|\bar{\kappa} \cos\theta - \bar{\tau} \sin\theta\right|, \; \tau = \bar{\kappa} \sin\theta + \bar{\tau} \cos\theta. \tag{73}$$

*(iv) If $\alpha$ is a spacelike curve of type 2, $X$ is spacelike and $\beta$ is a spacelike Bertrand-direction curve of type 2, then*

$$\bar{\kappa} = \left|\kappa \cosh\theta + \tau \sinh\theta\right|, \quad \bar{\tau} = -\kappa \sinh\theta - \tau \cosh\theta. \tag{74}$$

$$\kappa = \left|\bar{\kappa} \cosh\theta + \bar{\tau} \sinh\theta\right|, \; \tau = -\bar{\kappa} \sinh\theta - \bar{\tau} \cosh\theta. \tag{75}$$

*(v) If $\alpha$ is a spacelike curve of type 2, $X$ is timelike and $\beta$ is a timelike Bertrand-direction curve of $\alpha$, then*

$$\bar{\kappa} = \left|\kappa \sinh\theta + \tau \cosh\theta\right|, \quad \bar{\tau} = -\kappa \cosh\theta - \tau \sinh\theta. \tag{76}$$

$$\kappa = \left|-\bar{\kappa} \sinh\theta - \bar{\tau} \cosh\theta\right|, \tau = \bar{\kappa} \cosh\theta + \bar{\tau} \sinh\theta. \tag{77}$$

**Proof.** Let give the proof of case (i). From (52) and the second equation of system (57), we have

$$\bar{\kappa} = \left|\kappa \cosh\theta - \tau \sinh\theta\right|. \tag{78}$$

Moreover, from (67), we easily get

$$\bar{B}' = -\cosh\theta T + (-\kappa \sinh\theta + \tau \cosh\theta) N - \sinh\theta B.$$

Since $\bar{N} = N$ and $\alpha$ is timelike, we have

$$\bar{\tau} = \left\langle \bar{B}', \bar{N} \right\rangle = -\kappa \sinh\theta + \tau \cosh\theta. \tag{79}$$

From (78) and (79) we easily obtain (69).

The proofs of other cases can be given by a similar way of proof of case (i).

**Corollary 5.2.** *Let $\beta$ be non-lightlike Bertrand-direction curve of $\alpha$. Then one of the followings hold,*

*(i)* $\dfrac{\kappa^2}{\left(\kappa^2 - \tau^2\right)^{3/2}} \left(\dfrac{\tau}{\kappa}\right)' = \dfrac{\bar{\kappa}^2}{\left(\bar{\kappa}^2 - \bar{\tau}^2\right)^{3/2}} \left(\dfrac{\bar{\tau}}{\bar{\kappa}}\right)'$, *(ii)* $-\dfrac{\kappa^2}{\left(\kappa^2 - \tau^2\right)^{3/2}} \left(\dfrac{\tau}{\kappa}\right)' = \dfrac{\bar{\kappa}^2}{\left(\bar{\kappa}^2 - \bar{\tau}^2\right)^{3/2}} \left(\dfrac{\bar{\tau}}{\bar{\kappa}}\right)'$,

*(iii)* $\dfrac{\kappa^2}{\left(\kappa^2 + \tau^2\right)^{3/2}} \left(\dfrac{\tau}{\kappa}\right)' = \dfrac{\bar{\kappa}^2}{\left(\bar{\kappa}^2 + \bar{\tau}^2\right)^{3/2}} \left(\dfrac{\bar{\tau}}{\bar{\kappa}}\right)'$.

Now, we focus on relations between non-lightlike Bertrand-direction curves and some special curves such as helices and slant helices in $E_1^3$.

***Theorem 5.3.*** *Let $\alpha : I \to E_1^3$ be a non-lightlike curve and $\beta$ be a non-lightlike Bertrand-direction curve of $\alpha$. Then the following statements are equivalent*

*(i) $\alpha$ is a helix.*

*(ii) $\alpha$ is a Bertrand-donor curve of a helix.*

*(iii) A Bertrand-direction curve of $\alpha$ is a helix.*

***Theorem 5.4.*** *Let* $\alpha: I \to E_1^3$ *be a non-lightlike curve and* $\beta$ *be a non-lightlike Bertrand-direction curve of* $\alpha$. *If* $\alpha$ *is a plane curve, then* $\beta$ *is a helix. Similarly, if* $\beta$ *is a plane curve, then* $\alpha$ *is a helix.*

***Theorem 5.5.*** *Let* $\alpha: I \to E_1^3$ *be a non-lightlike curve and* $\beta$ *be a non-lightlike Bertrand-direction curve of* $\alpha$. *Then the following statements are equivalent*

*(i)* $\alpha$ *is a slant helix.*

*(ii)* $\alpha$ *is a Bertrand-donor curve of a slant helix.*

*(iii) A Bertrand-direction curve of* $\alpha$ *is a slant helix.*

**6. Example:** Consider the spacelike curve $\alpha$ of type 1 given by the parametric form $\alpha(s) = \frac{1}{2}\left(\cosh(s), \sinh(s), \sqrt{3}s\right)$. The Frenet vectors and curvatures of $\alpha$ are

$T(s) = \frac{1}{2}\left(\sinh(s), \cosh(s), \sqrt{3}\right)$, $N(s) = \left(\cosh(s), \sinh(s), 0\right)$, $B(s) = -\frac{1}{2}\left(\sqrt{3}\sinh(s), \sqrt{3}\cosh(s), 1\right)$

and $\kappa = 1/2,\ \tau = -\sqrt{3}/2$, respectively. It is clear that $\alpha$ is a general helix in $E_1^3$.

(i) From (14), an evolute-direction curve of $\alpha$ has the parametrization

$$\beta(s) = \int_0^s \beta'(s)ds = \int_0^s X(s)ds = \left(\beta_1(s), \beta_2(s), \beta_3(s)\right),$$

where

$$\beta_1(s) = \int_0^s \left[\sinh\left(-\frac{\sqrt{3}}{2}s + c\right)\cosh(s) + \frac{\sqrt{3}}{2}\cosh\left(-\frac{\sqrt{3}}{2}s + c\right)\sinh(s)\right]ds,$$

$$\beta_2(s) = \int_0^s \left[\sinh\left(-\frac{\sqrt{3}}{2}s + c\right)\sinh(s) + \frac{\sqrt{3}}{2}\cosh\left(-\frac{\sqrt{3}}{2}s + c\right)\cosh(s)\right]ds,$$

$$\beta_3(s) = \frac{1}{2}\int_0^s \cosh\left(-\frac{\sqrt{3}}{2}s + c\right)ds.$$

And, from Theorem 3.3, $\beta(s)$ is a spacelike slant helix.

(ii) From (28), a Mannheim-direction curve of $\alpha$ is obtained as follows

$$\gamma(s) = \int_0^s \gamma'(s)ds = \int_0^s X(s)ds = \left(\gamma_1(s), \gamma_2(s), \gamma_3(s)\right),$$

where

$$\gamma_1(s) = \int_0^s \left[\frac{1}{2}\sinh\left(\frac{s}{2} + c\right)\sinh(s) - \cosh\left(\frac{s}{2} + c\right)\cosh(s)\right]ds,$$

$$\gamma_2(s) = \int_0^s \left[\frac{1}{2}\sinh\left(\frac{s}{2} + c\right)\cosh(s) - \cosh\left(\frac{s}{2} + c\right)\sinh(s)\right]ds,$$

$$\gamma_3(s) = \frac{\sqrt{3}}{2}\int_0^s \sinh\left(\frac{s}{2} + c\right)ds.$$

From Theorem 4.3, $\gamma(s)$ is a timelike slant helix.

(iii) From (54), a Bertrand-direction curve of $\alpha$ is obtained as follows

$$\delta(s) = \int_0^s \delta'(s)ds = \int_0^s X(s)ds = \left(\delta_1(s), \delta_2(s), \delta_3(s)\right),$$

where

$$\delta_1(s)=\int_0^s\left[\frac{1}{2}\cos(\theta)\sinh(s)-\frac{\sqrt{3}}{2}\sin(\theta)\sinh(s)\right]ds,$$

$$\delta_2(s)=\int_0^s\left[\frac{1}{2}\cos(\theta)\cosh(s)-\frac{\sqrt{3}}{2}\sin(\theta)\cosh(s)\right]ds,$$

$$\delta_3(s)=\int_0^s\left[\frac{\sqrt{3}}{2}\cos(\theta)-\frac{1}{2}\sin(\theta)\right]ds,$$

where $\theta$ is constant and from Theorem 5.3, $\delta(s)$ is a spacelike helix.

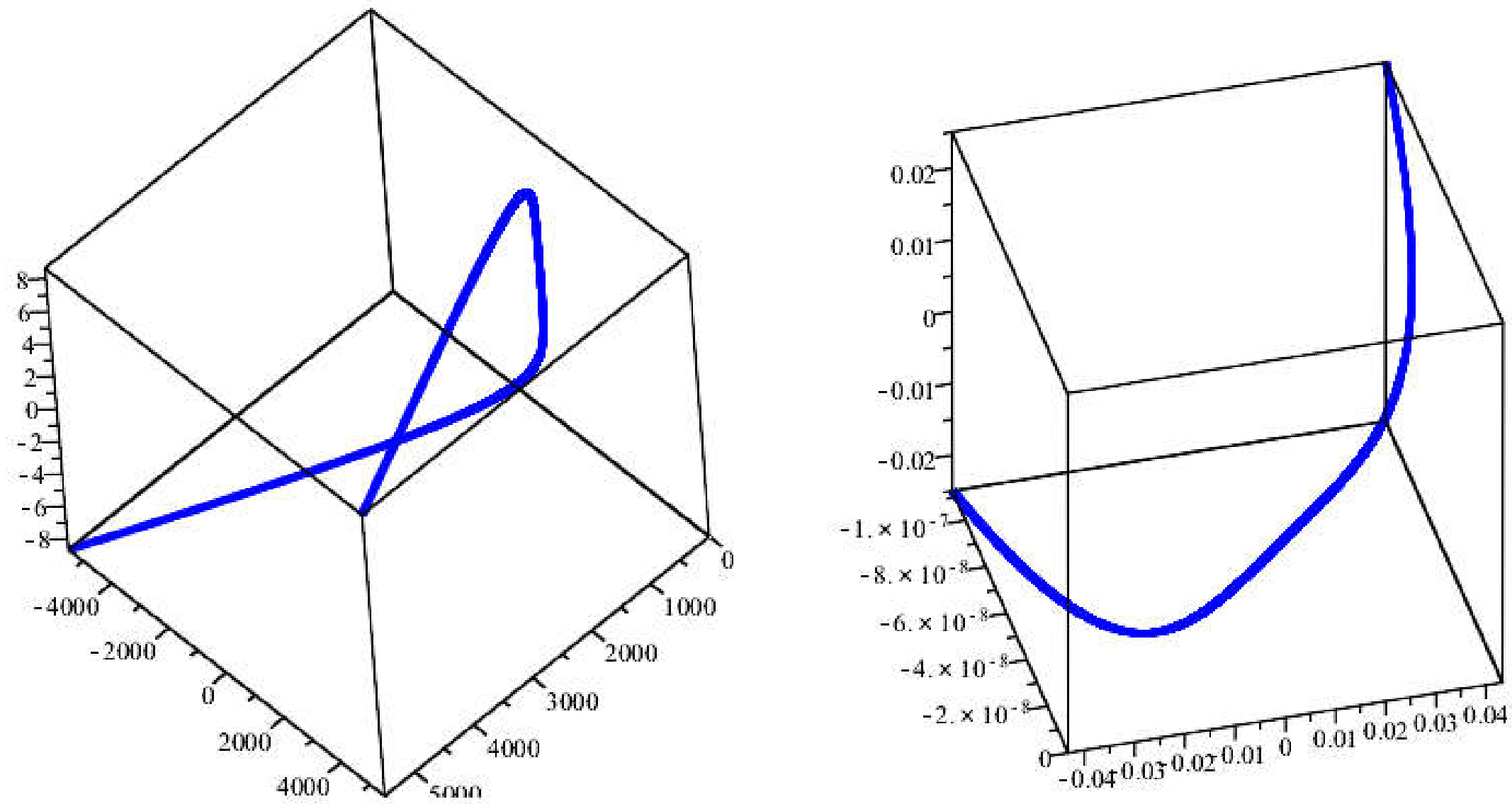

**Fig 1.** Spacelike helix $\alpha$ (left) and evolute-direction curve $\beta$ (right).

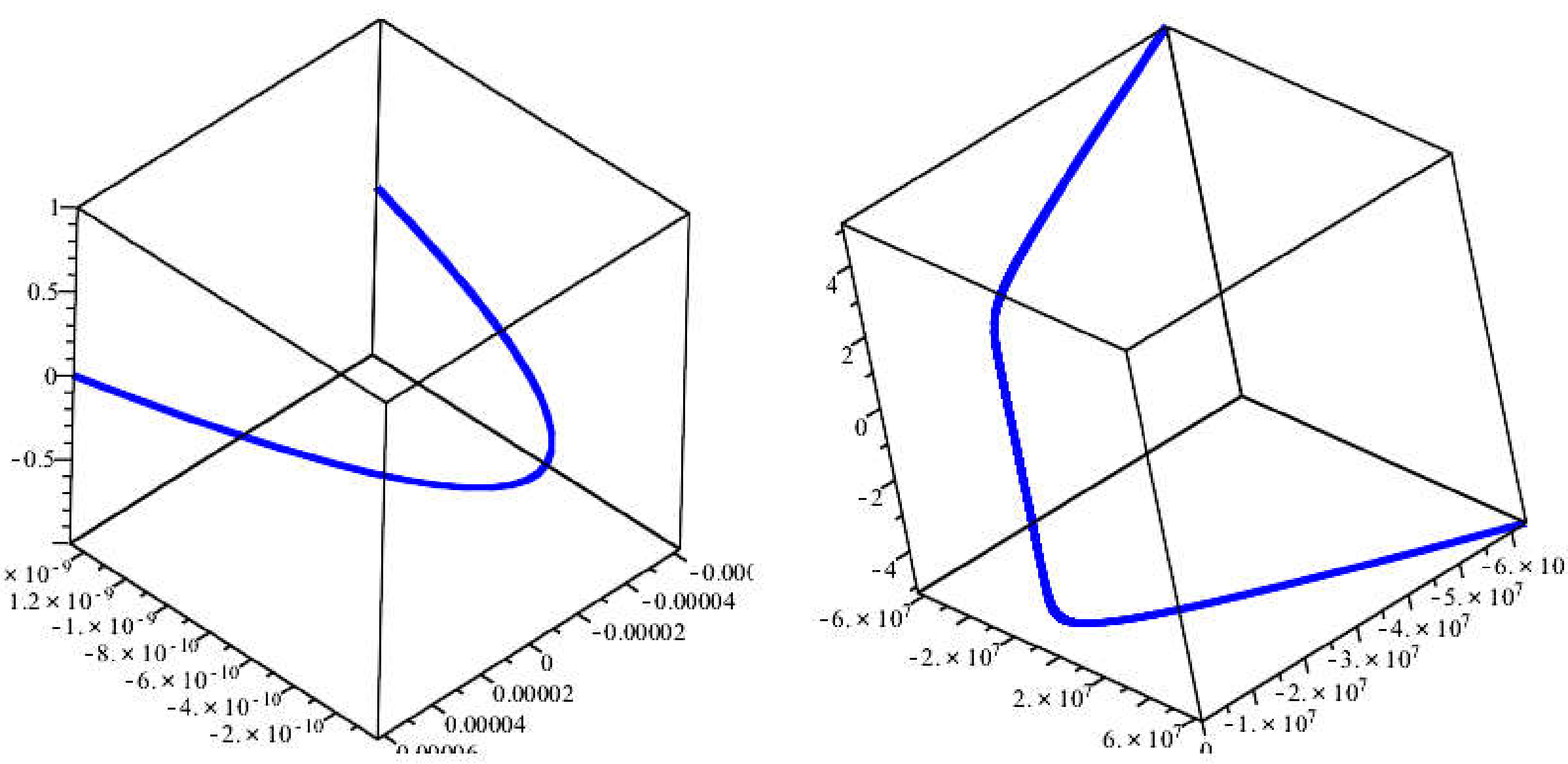

**Fig 2.** Mannheim-direction curve $\gamma$ (left) and Bertrand-direction curve $\delta$ (right).

## Conclusions

In the present paper, integral curves (an important subject of physics) are used to define and study some new curve couples such as involute-evolute-direction curves, Mannheim partner-direction curves and Bertrand partner-direction curves in Minkowski 3-space. The proved theorems give a way to obtain relations between physics and geometry. The most important insight and application is that integral curves (as a physical subject) give a way to construct some special curves such as helix and slant helix which are the most famous subject of differential geometry.